# Exponential functionals of Brownian motion, I: Probability laws at fixed time[*]

**Hiroyuki Matsumoto**

*Graduate School of Information Science, Nagoya University,
Chikusa-ku, Nagoya 464-8601, Japan
e-mail:* `matsu@info.human.nagoya-u.ac.jp`

**Marc Yor**

*Laboratoire de Probabilités and Institut universitaire de France,
Université Pierre et Marie Curie,
175 rue du Chevaleret, F-75013 Paris, France
e-mail:* `deaproba@proba.jussieu.fr`

**Abstract:** This paper is the first part of our survey on various results about the distribution of exponential type Brownian functionals defined as an integral over time of geometric Brownian motion. Several related topics are also mentioned.

**AMS 2000 subject classifications:** Primary 60J65.
**Keywords and phrases:** Brownian motion, Bessel process, Lamperti's relation, Hartman-Watson distributions.

## 1. Introduction

Throughout this paper we denote by $B = \{B_t, t \geqq 0\}$ a one-dimensional Brownian motion starting from 0 defined on a probability space $(\Omega, \mathcal{F}, P)$ and by $B^{(\mu)} = \{B_t^{(\mu)} = B_t + \mu t, t \geqq 0\}$ the corresponding Brownian motion with constant drift $\mu \in \mathbf{R}$. In this and the second part of our survey ([51]), we are concerned with the additive functional $A^{(\mu)} = \{A_t^{(\mu)}\}$ defined by

$$A_t^{(\mu)} = \int_0^t \exp(2B_s^{(\mu)})ds, \qquad t \geqq 0.$$

When $\mu = 0$, we simply write $A_t$ for $A_t^{(0)}$.

Such an exponential type functional plays an important role in several domains, e.g., mathematical finance, diffusion processes in random environments, stochastic analysis related to Brownian motions on hyperbolic spaces. This fact motivates detailed studies about this functional.

---

[*]This is an original survey paper.

In this first part, we mainly consider the distribution of $A_t^{(\mu)}$ at fixed time $t$, while some related results are also shown. It should be mentioned that, although $A_t^{(\mu)}$ is of a simple form, only integral representations for the density (of somewhat complicated forms!) are known by the results due to Alili-Gruet [3], Comtet-Monthus-Yor [15], Dufresne [23], Yor [69]. See also [8],[14], [44], [50], [52] and [61].

We mention the results obtained in [3], [15], [23] and [69] in some detail and we show some related identities, which are of independent interest. For example, the Bougerol identity $\sinh(B_t) \stackrel{\text{(law)}}{=} \beta_{A_t}$ for fixed $t$ and for another Brownian motion $\{\beta_u\}$ independent of $B$ and the Dufresne recursion relation, which relates the laws of $A_t^{(\mu)}$ and $A_t^{(\nu)}$ (different drifts), will be mentioned. We also show some results about the moments of $A_t^{(\mu)}$ and their asymptotics as $t \to \infty$.

The organization of this paper is as follows.

In Section 2 we introduce the Hartman-Watson distribution after recalling some facts on Bessel processes.

In Section 3 we show the Bougerol identity.

In Section 4 we present several expressions for the density of $A_t^{(\mu)}$ at a fixed time $t$ and at an independent exponential time.

In Section 5 we present the computation of the moments of $A_t^{(\mu)}$. One may derive or guess some important identities from the expressions for the moments, and we also present them.

In Section 6 we study various properties of perpetual functionals of the form $\int_0^\infty f(B_t^{(\mu)})dt$.

In Section 7, after giving some results on the asymptotic behavior of the Laplace transform and the moments of $A_t^{(\mu)}$ as $t$ tends to $\infty$, we show some limit theorems for some Gibbsian measures and moment density measures, following Hariya-Yor [32].

In Section 8, considering another exponential functional of Brownian motion, we present a result on a three-dimensional joint distribution.

Finally, in Appendix A, we give a proof of an integral representation for the density of the Hartman-Watson distribution.

## 2. Hartman-Watson distributions

In this section we give some information about Bessel processes and introduce the Hartman-Watson distributions, which were originally obtained and studied in [33] from an analytic point of view.

The starting important fact is the Lamperti relation, which says that there exists a Bessel process $R^{(\mu)} = \{R_u^{(\mu)}\}$ with $R_0^{(\mu)} = 1$ and index $\mu$ satisfying

$$\exp(B_t^{(\mu)}) = R_{A_t^{(\mu)}}^{(\mu)}, \quad t \geqq 0. \tag{2.1}$$

For details, see Lamperti [40] and also p.452 of [57].

When $\mu < 0$, we easily see that $\lim_{t\to\infty} A_t^{(\mu)} = A_\infty^{(\mu)}$ is finite and, from (2.1), that
$$A_\infty^{(\mu)} = \inf\{s; R_s^{(\mu)} = 0\} < \infty.$$

Dufresne [21] (see also [70]) proved that $A_\infty^{(\mu)}$ is distributed as $1/2\gamma_{-\mu}$, where $\gamma_{-\mu}$ is a gamma random variable with parameter $-\mu$. We will give a proof in Section 6 of this article. Baldi et al [4] have discussed about $A_\infty^{(\mu)}$ in relation to some diffusion processes on the hyperbolic plane.

In view of applications, it may also be worthwhile to note
$$A^{(\mu)}_{\tau_b(B^{(\mu)})} = \tau_{\exp(b)}(R^{(\mu)}),$$

where $\tau_c(X)$ denotes the first hitting time at $c \in \mathbf{R}$ of a stochastic process $X$. When $\mu > 0$, replacing $\tau_c$ by, e.g., last hitting times, we have a similar identity.

Now denote by $P_x^{(\mu)}$ the probability law of a Bessel process $R^{(\mu)}$ with index $\mu$ on the canonical path space $C([0,\infty) \to \mathbf{R})$ when $R_0^{(\mu)} = x$. Then, letting $\mathcal{R}_t = \sigma\{R_s, s \leqq t\}$ be the filtration of the coordinate process $\{R_u\}$, we have the following absolute continuity relationship for the laws of the Bessel processes with different indices: for $\mu \geqq 0$ and $x > 0$,
$$P_x^{(\mu)}\big|_{\mathcal{R}_t} = \left(\frac{R_t}{x}\right)^\mu \exp\left(-\frac{\mu^2}{2}\int_0^t \frac{ds}{(R_s)^2}\right) \cdot P_x^{(0)}\big|_{\mathcal{R}_t} \qquad (2.2)$$

and, for $\mu < 0$ and $x > 0$,
$$P_x^{(\mu)}\big|_{\mathcal{R}_t \cap \{t < \tau_0(R)\}} = \left(\frac{R_t}{x}\right)^\mu \exp\left(-\frac{\mu^2}{2}\int_0^t \frac{ds}{(R_s)^2}\right) \cdot P_x^{(0)}\big|_{\mathcal{R}_t}. \qquad (2.3)$$

We note that these absolute continuity relationships and an extension of formula (2.8) below also hold for the Wishart processes, which are diffusion processes with values in the space of positive definite matrices. For details, see [17].

The important relationships (2.2) and (2.3) are shown in the following way. Since formula (2.2) has been presented on p. 450, [57] and is proven by an easy modification of the arguments below, we only give a proof for (2.3) in the case $x = 1$.

Let $W^{(\mu)}, \mu \in \mathbf{R}$, be the probability law of $\{B_t^{(\mu)}\}$ on $C([0,\infty) \to \mathbf{R})$ and $\{\mathcal{F}_t\}$ be the filtration of the coordinate process $\{X_t\}$. Then, by the Cameron-Martin theorem, we have
$$W^{(\mu)}\big|_{\mathcal{F}_t} = \exp(2\mu X_t) W^{(-\mu)}\big|_{\mathcal{F}_t}.$$

We set $A_t(X) = \int_0^t \exp(2X_s)ds$, $\alpha_u = \inf\{t > 0; A_t(X) > u\}$ and define $R = \{R_u\}$ by $R_u = \exp(X_{\alpha_u})$. Then $\{R_u, u < \tau_0(R)\}$ is a Bessel process with index $\mu \in \mathbf{R}$. Note that $\{u < \tau_0(R)\} = \{\alpha_u < \infty\}$ and $0 < W^{(\mu)}(\alpha_u < \infty) < 1$ if $\mu < 0$. Therefore, by time change, we obtain
$$W^{(\mu)}\big|_{\mathcal{F}_{\alpha_u}} = (R_u)^{2\mu} W^{(-\mu)}\big|_{\mathcal{F}_{\alpha_u} \cap \{u < \tau_0(R)\}} \qquad \text{for } \mu > 0,$$

which, combined with (2.2), implies (2.3).

On the other hand, letting $p^{(\mu)}(t,x,y), t > 0, x, y \geqq 0$, be the transition probability density with respect to the Lebesgue measure of the Bessel process $R^{(\mu)}$, its explicit expression is also well known ([57], p.446): if $\mu > -1$,

$$p^{(\mu)}(t,x,y) = \frac{1}{t}\left(\frac{y}{x}\right)^\mu y \exp\left(-\frac{x^2+y^2}{2t}\right) I_\mu\left(\frac{xy}{t}\right) \qquad (2.4)$$

for $x > 0$ and, for $x = 0$,

$$p^{(\mu)}(t,0,y) = \frac{1}{2^\mu t^{1+\mu}\Gamma(1+\mu)} y^{1+2\mu} e^{-y^2/2t},$$

where $I_\mu$ is the modified Bessel function. For a proof of (2.4), see Remark 2.1 below.

From formula (2.4), assuming $R_0^{(\mu)} = x$, we obtain

$$E[\exp(-\lambda(R_t^{(\mu)})^2)] = \frac{1}{(1+2\lambda t)^{1+\mu}} \exp\left(-\frac{x^2\lambda}{1+2\lambda t}\right), \quad \lambda > 0. \qquad (2.5)$$

Note that, in fact, in [57], formula (2.4) has been deduced as a consequence of (2.5).

*Remark* 2.1. As usual, denote by $K_\mu$ the other modified Bessel (Macdonald) function. We set

$$F_\mu(x,y) = \begin{cases} I_\mu(x)K_\mu(y), & \text{if } x \leqq y, \\ K_\mu(x)I_\mu(y), & \text{if } x \geqq y. \end{cases} \qquad (2.6)$$

Then, by the general theory of the Sturm-Liouville operators, we can show that the Green function (with respect to the Lebesgue measure) is given by

$$\int_0^\infty e^{-\alpha t} p^{(\mu)}(t,x,y) dt = 2y\left(\frac{y}{x}\right)^\mu F_\mu(\sqrt{2\alpha}x, \sqrt{2\alpha}y)$$

and, combining this with the following formula for the product of the modified Bessel functions (cf. [25], p.53)

$$I_\mu(x)K_\mu(y) = \frac{1}{2}\int_0^\infty e^{-\xi/2-(x^2+y^2)/2\xi} I_\mu\left(\frac{xy}{\xi}\right) \frac{d\xi}{\xi}, \quad 0 < x \leqq y, \qquad (2.7)$$

we obtain (2.4).

The function $F_\mu$ will be used to give an expression for the double Laplace transform of the distribution of $A_t^{(\mu)}$ and, when we invert it, formula (2.7) will play a role. For details, see Section 4.

Comparing formulae (2.2) and (2.4), we obtain the important relationship

$$E_x^{(0)}\left[\exp\left(-\frac{\mu^2}{2}\int_0^t \frac{ds}{(R_s)^2}\right)\bigg|R_t = y\right] = \left(\frac{I_{|\mu|}}{I_0}\right)\left(\frac{xy}{t}\right), \qquad (2.8)$$

where $E_x^{(0)}$ denotes the expectation with respect to $P_x^{(0)}$.

The conditional distribution of $\int_0^t (R_s)^{-2} ds$ given $R_t = y$ under $P_x^{(0)}$, characterized by (2.8), is the Hartman-Watson distribution with parameter $r = xy/t$. Of course this is a probabilistic interpretation, and Hartman-Watson [33] have introduced a family of probability measures $\eta_r(dt)$, $r > 0$, on $\mathbf{R}_+$, from the identity

$$\int_0^\infty e^{-\mu^2 t/2} \eta_r(dt) = \left(\frac{I_{|\mu|}}{I_0}\right)(r), \quad \mu \in \mathbf{R}.$$

For more probabilistic discussions of the Hartman-Watson distributions, see Yor [67].

In [67] an integral representation for the density of $\eta_r(dt)$ has been obtained. Let us define the function $\theta(r,t)$, $r > 0$, $t > 0$, by

$$\theta(r,t) = \frac{r}{(2\pi^3 t)^{1/2}} e^{\pi^2/2t} \int_0^\infty e^{-\xi^2/2t} e^{-r\cosh(\xi)} \sinh(\xi) \sin\left(\frac{\pi\xi}{t}\right) d\xi. \quad (2.9)$$

Then it holds that $I_0(r)\eta_r(dt) = \theta(r,t) dt$, or equivalently

$$\int_0^\infty e^{-\mu^2 t/2} \theta(r,t) dt = I_{|\mu|}(r), \quad r > 0. \quad (2.10)$$

A proof of formula (2.10) will be given in Appendix A. We will see in Section 4 that an expression for the joint probability density of $(A_t^{(\mu)}, B_t^{(\mu)})$ for a fixed $t$ may be given by using the function $\theta(r,t)$.

From the explicit expression for $\theta(r,t)$, one deduces

$$\lim_{t \to \infty} \sqrt{2\pi t^3} \theta(r,t) = r \int_0^\infty \xi e^{-r\cosh(\xi)} \sinh(\xi) d\xi$$
$$= \int_0^\infty e^{-r\cosh(\xi)} d\xi = K_0(r). \quad (2.11)$$

*Remark* 2.2. Let $\{\varphi_t\}$ be a continuous determination of the angular component of a complex Brownian motion $\{Z_t\}$ with $Z_0 = z \neq 0$:

$$Z_t = |Z_t| \exp(\sqrt{-1}\varphi_t).$$

Then it holds ([67]) that

$$E[\exp(\sqrt{-1}\lambda(\varphi_t - \varphi_0)) \mid |Z_t| = \rho] = \left(\frac{I_{|\lambda|}}{I_0}\right)\left(\frac{|z|\rho}{t}\right). \quad (2.12)$$

Moreover it should also be noted that there exists a one-dimensional Brownian motion $\{\beta_u\}$ independent of $\{|Z_s|\}$ such that

$$\varphi_t - \varphi_0 = \beta_{\alpha_t} \quad \text{with} \quad \alpha_t = \int_0^t \frac{ds}{|Z_s|^2}. \quad (2.13)$$

See also Itô-McKean [37], p.270.

*Remark* 2.3. The function $\theta(r,t)$ is, of course, a non-negative function. However, since the integral on the right hand side of (2.9) is an oscillatory one, it is not easy to carry out numerical computations, especially when $t$ is small. See [5], [36] and [44] for some results and discussions on the numerical computations.

*Remark* 2.4. In his fundamental memoir [63], Stieltjes considered a similar integral to that on the right hand side of (2.9). He has shown that

$$\int_0^\infty e^{-\xi^2/2t} \sinh(n\xi) \sin\left(\frac{\pi\xi}{t}\right) d\xi = 0 \tag{2.14}$$

holds for any $n \in \mathbf{N}$ and for any $t > 0$. Setting

$$\mu_{t,\lambda}(d\xi) = \frac{1}{\sqrt{2\pi t}} e^{-\xi^2/2t}\left(1 + \lambda \sin\left(\frac{\pi\xi}{t}\right)\right) d\xi, \qquad |\lambda| < 1,$$

we see from (2.14) that the exponential moments of $\mu_{t,\lambda}(d\xi)$ are given by

$$\int_{\mathbf{R}} e^{n\xi} \mu_{t,\lambda}(d\xi) = e^{n^2 t/2}, \qquad n = 1, 2, ...,$$

and that the moment problem for the log-normal distribution is indeterminate.

## 3. Bougerol's identity

In his discussion of stochastic analysis on hyperbolic spaces, Bougerol [11] has obtained an interesting and important identity in law. It plays an important role in several domains and in the following discussions. We show it in the simplest form. Some extensions and discussions are found in [2] and [3].

**Theorem 3.1.** *Let $\{W_t\}$ be a one-dimensional Brownian motion starting from $0$, independent of $B$. Then, for every fixed $t > 0$, $\int_0^t \exp(B_s)dW_s$, $\sinh(B_t)$ and $W_{A_t}$ are identical in law.*

*Proof.* We follow the proof in [2]. We put

$$X_t^x = \exp(B_t)\left(x + \int_0^t \exp(-B_s)dW_s\right)$$

and

$$Y_t^y = \sinh(y + B_t).$$

Then $\{X_t^x\}$ and $\{Y_t^y\}$ are diffusion processes with the same generator

$$\frac{1}{2}(1+x^2)\frac{d^2}{dx^2} + \frac{1}{2}x\frac{d}{dx}$$

and, if $\sinh(y) = x$, $\{X_t^x\}$ and $\{Y_t^y\}$ have the same probability law. Hence, we obtain the result if we consider the case $x = y = 0$ and if we notice that, for a fixed $t$, $X_t^0$ is identical in law with $\int_0^t \exp(B_s)dW_s$ by virtue of the invariance of the probability law of Brownian motion under time reversal from a fixed time. □

From Bougerol's identity $\sinh(B_t) \stackrel{(\text{law})}{=} W_{A_t}$, we easily obtain the following by considering the corresponding densities and characteristic functions.

**Corollary 3.2.** *For every $\alpha \in \mathbf{R}$ and $t > 0$, one has*

$$E\left[\frac{1}{\sqrt{A_t}}\exp\left(-\frac{\alpha^2}{2A_t}\right)\right] = \frac{1}{\sqrt{(1+\alpha^2)t}}\exp\left(-\frac{(\operatorname{Argsinh}(\alpha))^2}{2t}\right) \quad (3.1)$$

*and*

$$E\left[\exp\left(-\frac{\alpha^2}{2}A_t\right)\right] = \sqrt{\frac{2}{\pi t}}\int_0^\infty \cos(\alpha\sinh(\xi))e^{-\xi^2/2t}d\xi. \quad (3.2)$$

## 4. The law of $A_t^{(\mu)}$ at fixed and at exponential times

Since the exponential functional appears in several domains and the expressions for the density obtained so far are complicated, some authors have tried to obtain simpler expressions. We give a survey on the results due to Alili-Gruet [3], Comtet-Monthus-Yor [15], Dufresne [23] and Yor [69]. For other approaches, see Bhattacharya-Thomann-Waymire [8], Comtet-Monthus [14], Lyasoff [44], Monthus [52], Schröder [61] and so on.

**4.1.** The first result on an explicit expression for the probability density of $A_t^{(\mu)}$ was given by Yor [69].

**Theorem 4.1.** *Fix $t > 0$. Then, for $u > 0$ and $x \in \mathbf{R}$, it holds that*

$$P(A_t^{(\mu)} \in du, B_t^{(\mu)} \in dx) = e^{\mu x - \mu^2 t/2}\exp\left(-\frac{1+e^{2x}}{2u}\right)\theta(e^x/u, t)\frac{dudx}{u}, \quad (4.1)$$

*where $\theta$ is the function given by (2.9).*

*Proof.* We sketch a proof with $\mu = 0$. The general case is easily derived from the Cameron-Martin theorem once we have shown this special case. Moreover, since a probabilistic proof is given in [69], which is also found in [66], we have chosen to give a proof based on the general theory of the Sturm-Liouville operators.

For this purpose, let us consider the Schrödinger operator $H_\lambda$, $\lambda > 0$, on $\mathbf{R}$ with the so-called Liouville potential given by

$$H_\lambda = -\frac{1}{2}\frac{d^2}{dx^2} + \frac{1}{2}\lambda^2 e^{2x}, \qquad x \in \mathbf{R}. \quad (4.2)$$

Then, by using some properties of the modified Bessel functions, we can show that the Green function $G(x, y; \alpha^2/2) = (H_\lambda + \alpha^2/2)^{-1}(x, y)$, $\alpha > 0$, with respect to the Lebesgue measure is given by

$$G(x, y; \alpha^2/2) = 2I_\alpha(\lambda e^x)K_\alpha(\lambda e^y) = 2F_\alpha(\lambda e^x, \lambda e^y), \quad x \leqq y, \quad (4.3)$$

where the function $F_\alpha$ is given by (2.6).

Next we let $q_\lambda(t, x, y)$ be the heat kernel of the semigroup $\exp(-tH_\lambda)$, $t > 0$. Then we have by the Feynman-Kac formula

$$q_\lambda(t, x, y) = E\left[\exp\left(-\frac{1}{2}\lambda^2 e^{2x} A_t\right)\bigg| B_t = y - x\right] \frac{1}{\sqrt{2\pi t}} e^{-(y-x)^2/2t}$$
$$= \int_0^\infty e^{-\lambda^2 e^{2x} u/2} P(A_t \in du | B_t = y - x) \frac{1}{\sqrt{2\pi t}} e^{-(y-x)^2/2t} \quad (4.4)$$

and

$$\int_0^\infty e^{-\alpha^2 t/2} q_\lambda(t, x, y) dt = G(x, y; \alpha^2/2).$$

Now recall the integral representation (2.7) for the product of two modified Bessel functions. Then, combining the formulae above and setting $x = 0$, we obtain from (2.10)

$$\int_0^\infty \frac{1}{\sqrt{2\pi t}} e^{-y^2/2t - \alpha^2 t/2} dt \int_0^\infty e^{-\lambda^2 u/2} P(A_t \in du | B_t = y)$$
$$= \int_0^\infty \exp\left(-\frac{1}{2}u - \frac{\lambda^2(1 + e^{2y})}{2u}\right) I_\alpha\left(\frac{\lambda^2 e^y}{u}\right) \frac{du}{u}$$
$$= \int_0^\infty \exp\left(-\frac{1}{2}u - \frac{\lambda^2(1 + e^{2y})}{2u}\right) \frac{du}{u} \int_0^\infty e^{-\alpha^2 t/2} \theta\left(\frac{\lambda^2 e^y}{u}, t\right) dt.$$

Using the Fubini theorem, we can easily invert the double Laplace transform and obtain

$$P(A_t \in du | B_t = y) \frac{1}{\sqrt{2\pi t}} e^{-y^2/2t} = \exp\left(-\frac{1 + e^{2y}}{2u}\right) \theta(e^y/u, t) \frac{du}{u},$$

which is equivalent to (4.1) with $\mu = 0$. □

*Remark* 4.1. It may be worth presenting an explicit expression for the heat kernel $q_\lambda(t, x, y)$, which is obtained from (4.4):

$$q_\lambda(t, x, y) = \int_0^\infty \exp\left(-\frac{\xi}{2} - \frac{\lambda^2(e^{2x} + e^{2y})}{2\xi}\right) \theta\left(\frac{\lambda^2 e^{x+y}}{\xi}, t\right) \frac{d\xi}{\xi}.$$

Another expression will be given by Proposition 5.4 below.

*Remark* 4.2. We find formula (4.3) in Borodin-Salminen [9], page 264. In that book, many explicit formulae related to the exponential functionals are found. Moreover, some physicists working with Feynman path integrals have obtained similar results in the context of physical problems. We find many formulae and related references in Grosche-Steiner [28]. For example, on page 228 of [28], we find a formula equivalent to (4.3).

We obtain an expression for the density of $A_t^{(\mu)}$ from (4.1). However, to compare with Dufresne's results mentioned in Subsection 4.3, it is convenient to consider $1/2A_t^{(\mu)}$. We denote its density by $f^{(\mu)}(a,t)$. Then, by using the integral representation (2.9) for the function $\theta(r,t)$, we have

$$f^{(\mu)}(a,t) = \frac{2}{(2\pi^3 t)^{1/2}} e^{\pi^2/2t - \mu^2 t/2} e^{-a} a^{-(\mu+1)/2} \int_0^\infty \eta^\mu e^{-\eta^2} d\eta$$
$$\times \int_0^\infty e^{-\xi^2/2t} e^{-2\sqrt{a}\eta\cosh(\xi)} \sinh(\xi) \sin\left(\frac{\pi\xi}{t}\right) d\xi$$

if $\mu > -1$.

When $\mu = 0$, using the elementary formula

$$e^{-\sqrt{\lambda}\zeta} = \frac{\zeta}{\sqrt{2}} \int_0^\infty \exp\left(-\lambda s - \frac{\zeta^2}{4s}\right) \frac{ds}{\sqrt{2\pi s^3}}, \quad \lambda, \zeta > 0,$$

we obtain

$$f^{(0)}(a,t) = \frac{2e^{\pi^2/2t}}{\pi^2 (2t)^{1/2}} \frac{1}{\sqrt{a}} \int_0^\infty e^{-a(\cosh(\eta))^2} \cosh(\eta) d\eta$$
$$\times \int_0^\infty e^{-\xi^2/2t} \frac{\sinh(2\xi)\sin(\pi\xi/t)}{\cosh(2\xi) + \cosh(2\eta)} d\xi. \quad (4.5)$$

There is also a similar expression for $f^{(\mu)}(a,t)$ when $\mu = 1$ (see [50]).

Moreover, by using the formula (4.1), we can compute the Laplace transform of the function $\theta(r,t)$ in $r$.

**Proposition 4.2.** Set $a_c(x) = \mathrm{Argcosh}(x)$ for $x \geq 1$. Then it holds that

$$\int_0^\infty e^{-xr} \theta(r,t) \frac{dr}{r} = \frac{1}{\sqrt{2\pi t}} \exp\left(-\frac{a_c(x)^2}{2t}\right) \quad (4.6)$$

for any $t > 0$. Consequently, one has

$$\int_0^\infty e^{-xr} \theta(r,t) dr = \frac{1}{\sqrt{2\pi t^3}} a_c(x) a_c'(x) \exp\left(-\frac{a_c(x)^2}{2t}\right). \quad (4.7)$$

*Proof.* Integrating both hand sides of (4.1) in $u$ yields

$$\frac{1}{\sqrt{2\pi t}} e^{-x^2/2t} = \int_0^\infty \exp\left(-\frac{1+e^{2x}}{2u}\right) \theta(e^x/u, t) \frac{du}{u}$$
$$= \int_0^\infty e^{-r\cosh(x)} \theta(r,t) \frac{dr}{r},$$

which is exactly (4.6).

Formula (4.7) follows from (4.6) by differentiating with respect to $x$. □

One finds another proof of Proposition 4.2 in [47], where the Lipschitz-Hankel formula below (cf. Chapter XIII, Watson [64]) has been used:

$$\int_0^\infty e^{-r\cosh(u)} I_\lambda(r) dr = \frac{1}{\sinh(u)} e^{-\lambda u}, \quad u > 0, \lambda > -1.$$

Several subordinators whose Lévy measures are written in terms of the modified Bessel functions have been studied in [47] in relation to the exponential functional of Brownian bridge.

Moreover, formula (4.6) may be proven from the explicit expression (2.9) for the function $\theta(r,t)$. In fact, by Fubini's theorem, we obtain

$$\int_0^\infty e^{-xr}\theta(r,t)\frac{dr}{r} = \frac{e^{\pi^2/2t}}{2\sqrt{2\pi^3 t}}\int_{\mathbf{R}} \frac{e^{-\xi^2/2t}\sinh(\xi)\sin(\pi\xi/t)}{x+\cosh(\xi)}d\xi$$

and, applying residue calculus to the integral on the right hand side, we obtain (4.6). By this method, we also obtain for $-1 < x \leq 1$

$$\int_0^\infty e^{-xr}\theta(r,t)\frac{dr}{r} = \frac{1}{\sqrt{2\pi t}}\exp\left(\frac{(\text{Arccos}(x))^2}{2t}\right).$$

**4.2.** We computed the double Laplace transform of the conditional distribution of $A_t^{(\mu)}$ under the condition $B_t^{(\mu)} = x$ and obtained an expression for the joint density of $(A_t^{(\mu)}, B_t^{(\mu)})$ in the previous subsection. Alili-Gruet [3] have obtained an explicit form of the double Laplace transform of the distribution of $A_t^{(\mu)}$.

**Theorem 4.3.** *For every $\mu \in \mathbf{R}$, one has*

$$\int_0^\infty e^{-\alpha^2 t/2} E\left[\exp\left(-\frac{1}{2}u^2 A_t^{(\mu)}\right)\right] dt = 2\int_{-\infty}^\infty e^{\mu y} F_\lambda(u, ue^y) dy$$

*for $\alpha, u > 0$, where $\lambda = \sqrt{\alpha^2 + \mu^2}$ and $F_\lambda$ is the function given by (2.6).*

*Proof.* We use the facts mentioned in Remarks 2.1 and 2.2. At first, by using the Cameron-Martin theorem, we rewrite the left hand side as follows:

$$\int_0^\infty e^{-\alpha^2 t/2} E\left[\exp\left(-\frac{1}{2}u^2 A_t^{(\mu)}\right)\right] dt$$
$$= \int_0^\infty e^{-(\alpha^2+\mu^2)t/2} E\left[\exp\left(\mu B_t - \frac{1}{2}u^2 A_t\right)\right] dt$$
$$= \int_0^\infty E\left[\exp\left(\mu B_t - \frac{1}{2}u^2 A_t + \sqrt{-1}\lambda\beta_t\right)\right] dt$$
$$= E\left[\int_0^\infty \exp\left(\mu B_t - \frac{1}{2}u^2 A_t + \sqrt{-1}\lambda\beta_t\right) dt\right],$$

where $\{\beta_t\}$ is a Brownian motion independent of $B$.

We recall the Lamperti relation and make the change of variables by $A_t = v$. Then, using the notations in Remark 2.2 and noting that $\{\alpha_v\}$ is the inverse process of $\{A_t\}$, we obtain

$$E\left[\int_0^\infty \exp\left(\mu B_t - \frac{1}{2}u^2 A_t + \sqrt{-1}\lambda\beta_t\right) dt\right]$$
$$= E\left[\int_0^\infty (R_v)^\mu \exp\left(-\frac{u^2 v}{2} + \sqrt{-1}\lambda\varphi_v\right) \frac{dv}{(R_v)^2}\right],$$

where $\{R_t\}$ is the radial process of a complex Brownian motion $Z$ starting from 1 and $\varphi_0 = 0$. Now, applying formula (2.12), we obtain

$$\int_0^\infty e^{-\alpha^2 t/2} E\left[\exp\left(-\frac{1}{2}u^2 A_t^{(\mu)}\right)\right] dt$$
$$= \int_0^\infty e^{-u^2 v/2} E[(R_v)^{\mu-2} E[e^{\sqrt{-1}\lambda\varphi_v}|\sigma\{R_v\}]] dv$$
$$= \int_0^\infty e^{-u^2 v/2} E\left[(R_v)^{\mu-2}\left(\frac{I_\lambda}{I_0}\right)\left(\frac{R_v}{v}\right)\right] dv.$$

Finally, using formulae (2.4) and (2.7), we obtain the theorem. □

Monthus [52] has shown another expression for the Laplace transform of the probability law of $A_t^{(\mu)}$. See also [3] for further discussions.

Moreover, Alili-Gruet have applied the Plancherel identity to the right hand side of formula (3.2) to obtain

$$E\left[\exp\left(-\frac{\alpha^2}{2}A_t\right)\right] = \frac{2}{\pi}\int_0^\infty e^{-\eta^2 t/2} \cosh\left(\frac{\pi}{2}\eta\right) K_{\sqrt{-1}\eta}(\alpha) d\eta,$$

which has been shown by other methods in [53] and [55].

**Theorem 4.4.** *One has*

$$P(A_t \in du) = \frac{du}{\sqrt{2\pi u^3}} \frac{1}{\sqrt{2\pi t}} \int_\mathbf{R} \cosh(\xi) e^{-(\cosh(\xi))^2/2u - (\xi + \sqrt{-1}\pi/2)^2/2t} d\xi$$

*for any $t > 0$ and*

$$\lim_{t \to \infty} \sqrt{2\pi t} P(A_t \in du) = \frac{1}{u} e^{-1/2u} du, \qquad u > 0.$$

To show the theorem, Alili-Gruet [3] have used the Plancherel identity and Comtet-Monthus-Yor [15] have used an integral representation for $K_{\sqrt{-1}\eta}$. For details, see the original papers.

**4.3.** Dufresne [23] has studied in detail the function

$$h^{\mu,r}(s,t) = e^{\mu^2 t/2} E\left[(2A_t^{(\mu)})^{-r} \exp\left(\frac{s}{2A_t^{(\mu)}}\right)\right], \tag{4.8}$$

which has been shown to be convergent for $r \in \mathbf{R}, \mathrm{Re}(s) < 1$ and $t > 0$. He has obtained a differential equation in $(s,t)$ which is satisfied by this function, explicit expressions in some special cases and a recurrence relation for $h^{\mu,r}(s,t)$ in $\mu$. From these results, he has obtained explicit expressions for the density of $A_t^{(\mu)}$.

At first we show a computation for the Laplace transform of $h^{\mu,r}(s,t)$ in $t$ by starting from formula (4.18) below, the notation of which we borrow. While this way is converse to Dufresne's original arguments, it is easier to understand. For $\text{Re}(s) \leq 0$ and $\lambda > \mu^2/2$, one has

$$\begin{aligned}
\bar{h}^{\mu,r}(s,\lambda) &\equiv \int_0^\infty e^{-\lambda t} h^{\mu,r}(s,t) dt \\
&= \int_0^\infty e^{-(\lambda-\mu^2/2)t} E\left[(2A_t^{(\mu)})^{-r} \exp\left(\frac{s}{2A_t^{(\mu)}}\right)\right] dt \\
&= \frac{1}{\lambda - \mu^2/2} E\left[\left(\frac{Z_\beta}{Z_{1,\alpha}}\right)^r \exp\left(\frac{sZ_\beta}{Z_{1,\alpha}}\right)\right] \\
&= \frac{1}{\lambda - \mu^2/2} \frac{1}{\Gamma(\beta)} \int_0^1 \alpha(1-x)^{\alpha-1} dx \int_0^\infty \left(\frac{y}{x}\right)^r e^{sy/x} y^{\beta-1} e^{-y} dy,
\end{aligned}$$

where $\alpha = (\sqrt{2\lambda}+\mu)/2$ and $\beta = (\sqrt{2\lambda}-\mu)/2$. After carrying out the integration in $y$ and making the change of variables by $x = 1 - v$, we can use an integral representation for the Gauss hypergeometric function $_2F_1$ (see, e.g., [41], p.240)

$$_2F_1(\alpha,\beta;\gamma;z) = \frac{\Gamma(\gamma)}{\Gamma(\beta)\Gamma(\gamma-\beta)} \int_0^1 t^{\beta-1}(1-t)^{\gamma-\beta-1}(1-tz)^{-\alpha} dt, \; \gamma > \beta > 0,$$

and we obtain

$$\bar{h}^{\mu,r}(s,\lambda) = \frac{\Gamma(\alpha)\Gamma(\beta+r)}{2(1-s)^{\beta+r}\Gamma(\alpha+\beta+1)} \; _2F_1\left(\alpha,\beta+r;\alpha+\beta+1;\frac{1}{1-s}\right). \quad (4.9)$$

After a careful inspection, Dufresne has shown this result under a weaker condition on the parameters.

**Theorem 4.5.** *If $\sqrt{2\lambda} > \max(-\mu, \mu - 2r)$ and if either* (i) $\text{Re}(s) < 0, r \in \mathbf{R}$, *or* (ii) $\text{Re}(s) = 0, r < 1$, *then formula (4.9) holds.*

In some special cases, the right hand side of (4.9) may be written by using elementary functions and we can invert the Laplace transform.

**Corollary 4.6.** *When $r = 1/2$ and $\mu = 0, 1$, one has, for $\text{Re}(s) \leq 0$,*

$$h^{0,1/2}(s,t) = \frac{1}{\sqrt{2t(1-s)}} \exp\left(-\frac{(\text{Argsinh}\sqrt{-s})^2}{2t}\right) \quad (4.10)$$

*and*

$$h^{1,1/2}(s,t) = \frac{1}{\sqrt{2t}} \exp\left(-\frac{(\text{Argsinh}\sqrt{-s})^2}{2t}\right). \quad (4.11)$$

*Remark* 4.3. Formula (4.10) is exactly the same as (3.1) when $s = -\alpha^2$.

We now try to invert the Laplace transform (4.8) of the probability law of $(A_t^{(\mu)})^{-1}$. Denote by $f^{(\mu)}(a,t)$ the density of $(2A_t^{(\mu)})^{-1}$. Then the inversion formula is

$$f^{(\mu)}(a,t) = e^{-\mu^2 t/2} \frac{a^{-r}}{2\pi\sqrt{-1}} \int_C e^{sx} h^{\mu,r}(-s,t) ds,$$

where the path of integration $C$ from $c - \sqrt{-1}\infty$ to $c + \sqrt{-1}\infty$ with $c \in \mathbf{R}$ is on the right side of the singularities of $h^{\mu,r}(-s,t)$. Note that the left hand side does not depend on $r$, while the right hand side apparently does.

When $\mu = 0$ and $\mu = 1$, Dufresne [23] has succeeded in carrying out the complex integral and has obtained the following fairly simple expressions for $f^{(\mu)}(a,t)$.

**Proposition 4.7.** *For $t > 0$ and $a > 0$, one has*

$$f^{(0)}(a,t) = \frac{2e^{\pi^2/8t}}{\pi\sqrt{2ta}} \int_0^\infty e^{-y^2/2t} e^{-a(\cosh(y))^2} \cosh(y) \cos\left(\frac{\pi y}{2t}\right) dy \quad (4.12)$$

*and*

$$f^{(1)}(a,t) = \frac{2e^{\pi^2/8t - t/2}}{\pi\sqrt{2ta}} \int_0^\infty e^{-y^2/2t} e^{-a(\cosh(y))^2} \sinh(y)\cosh(y) \sin\left(\frac{\pi y}{2t}\right) dy. \quad (4.13)$$

When $\mu = 0$, we can show the coincidence of the two expressions (4.5) and (4.12) for $f^{(0)}(a,t)$ by applying residue calculus to the right hand side of (4.5). We can also do the same thing when $\mu = 1$. Combining this coincidence with the recurrence formula (4.14) below also due to Dufresne [23] or (4.15), we obtain the desired coincidence for general values of $\mu$.

**Theorem 4.8.** *For all $\mu, r \in \mathbf{R}$, $s$ with $\mathrm{Re}(s) < 1$ and $t > 0$, one has*

$$h^{\mu,r}(s,t) = (1-s)^{\mu-r} h^{2r-\mu,r}(s,t). \quad (4.14)$$

**Theorem 4.9.** *For $\mu \in \mathbf{R}, a > 0$ and $t > 0$, one has*

$$f^{(\mu)}(a,t) = \frac{e^{\pi^2/8t - \mu^2 t/2}}{\pi(2a)^{(\mu+1)/2}\sqrt{t}} \int_{-\infty}^\infty e^{-y^2/2t} e^{-a(\cosh(y))^2} \cosh(y)$$
$$\times \cos\left(\frac{\pi}{2}\left(\frac{y}{t} - \mu\right)\right) H_\mu(\sqrt{a}\sinh(y)) dy,$$

*where $H_\mu$ is the Hermite function with parameter $\mu$.*

*Remark* 4.4. We refer the reader to Lebedev [41] for the Hermite functions. We do not need the assertion of Theorem 4.9 in the following and we skip explanations. For details of proofs of the theorems above, see the original paper [23].

Dufresne's recursion relation (4.14) has been extended to the level of stochastic processes and has been written in a more symmetric manner in [50]. The following theorem plays an important role in studying various diffusion processes related to the exponential functional $\{A_t^{(\mu)}\}$. For details and a proof, see [50] and [51].

**Theorem 4.10.** *Let $\mu < \nu$ and consider two Brownian motions $B^{(\mu)}, B^{(\nu)}$ with different drifts and the corresponding exponential additive functionals $\{A_t^{(\mu)}\}$,* $\{A_t^{(\nu)}\}$. Then, for every $t > 0$ and for every non-negative functional $F$ on $C((0, t] \to \mathbf{R})$, it holds that

$$e^{\mu^2 t/2} E\left[F\left(\frac{1}{A_s^{(\mu)}}, s \leqq t\right)\left(\frac{1}{A_t^{(\mu)}}\right)^m\right]$$
$$= e^{\nu^2 t/2} E\left[F\left(\frac{1}{A_s^{(\nu)}} + 2\gamma_\delta, s \leqq t\right)\left(\frac{1}{A_t^{(\nu)}}\right)^m\right],$$

where $m = (\mu + \nu)/2, \delta = (\nu - \mu)/2$ and $\gamma_\delta$ is a gamma random variable independent of $B^{(\nu)}$. In particular, for every non-negative function $\psi$ on $\mathbf{R}_+$, one has

$$e^{\mu^2 t/2} E\left[\psi\left(\frac{1}{A_t^{(\mu)}}\right)\left(\frac{1}{A_t^{(\mu)}}\right)^m\right] = e^{\nu^2 t/2} E\left[\psi\left(\frac{1}{A_t^{(\nu)}} + 2\gamma_\delta\right)\left(\frac{1}{A_t^{(\nu)}}\right)^m\right]$$

and

$$e^{\mu^2 t/2} f^{(\mu)}(a, t) = e^{\nu^2 t/2} \frac{1}{\Gamma(\delta)} a^{-m} e^{-a} \int_0^a (a-b)^{\delta-1} b^m e^b f^{(\nu)}(b, t) db. \quad (4.15)$$

**4.4.** Let $T_\lambda$ be an exponential random variable with parameter $\lambda > 0$ independent of the Brownian motion $B$. In [68], [69] and some other articles contained in [66], several identities for the distributions of $(A_{T_\lambda}^{(\mu)}, B_{T_\lambda}^{(\mu)})$ and $A_{T_\lambda}^{(\mu)}$ itself have been shown. Among them we present two results.

At first we present an explicit expression for the joint probability density of $(A_{T_\lambda}^{(\mu)}, B_{T_\lambda}^{(\mu)})$.

**Theorem 4.11.** *For $\mu \in \mathbf{R}$ and $\lambda > 0$, set $\nu = \sqrt{2\lambda + \mu^2}$ and let $p^{(\nu)}(t, x, y)$ be the transition probability density with respect to the Lebesgue measure of a Bessel process with index $\nu$. Then one has*

$$P(\exp(B_{T_\lambda}^{(\mu)}) \in dy, A_{T_\lambda}^{(\mu)} \in du) = \frac{\lambda}{y^{2+\nu-\mu}} p^{(\nu)}(u, 1, y) dy du. \quad (4.16)$$

*Proof.* At first we consider the case $\mu = 0$. Letting $f, g : \mathbf{R}_+ \to \mathbf{R}_+$ be non-negative Borel functions, we have

$$\Phi_\lambda^{(0)}(f, g) \equiv E[f(\exp(B_{T_\lambda})) g(A_{T_\lambda})]$$
$$= \lambda E\left[\int_0^\infty e^{-\lambda t} f(\exp(B_t)) g(A_t) dt\right]. \quad (4.17)$$

Recall the Lamperti relation (2.1) and, as in Subsection 4.2, put

$$\alpha_s = \int_0^s \frac{du}{(R_u^{(0)})^2},$$

which is the inverse process of $A_t$. Then, changing the variable in (4.17) by $t = \alpha_s$, we obtain

$$\Phi_\lambda^{(0)}(f, g) = \lambda E\left[\int_0^\infty e^{-\lambda H_s} f(R_s^{(0)}) g(s) \frac{ds}{(R_s^{(0)})^2}\right].$$

Moreover, by using (2.8) and then (2.4), we obtain

$$\Phi_\lambda^{(0)}(f,g) = \lambda \int_0^\infty E\left[\left(\frac{I_{\sqrt{2\lambda}}}{I_0}\right)\left(\frac{R_s^{(0)}}{s}\right)f(R_s^{(0)})\frac{1}{(R_s^{(0)})^2}\right]g(s)ds$$

$$= \lambda \int_0^\infty g(s)ds \int_0^\infty \left(\frac{I_{\sqrt{2\lambda}}}{I_0}\right)\left(\frac{y}{s}\right)f(y)\frac{1}{y^2}p^{(0)}(s,1,y)dy$$

$$= \lambda \int_0^\infty g(s)ds \int_0^\infty f(y)y^{-(2+\sqrt{2\lambda})}p^{(\sqrt{2\lambda})}(s,1,y)dy,$$

which implies formula (4.16) when $\mu = 0$.

Next we consider the general case. By the Cameron-Martin theorem, we have

$$\Phi_\lambda^{(\mu)}(f,g) \equiv E[f(\exp(B_{T_\lambda}^{(\mu)}))g(A_{T_\lambda}^{(\mu)})]$$

$$= \lambda \int_0^\infty e^{-\lambda t} E[f(\exp(B_t))g(A_t)\exp(\mu B_t - \mu^2 t/2)]dt$$

$$= \lambda \int_0^\infty e^{-\nu^2 t/2} E[f(\exp(B_t))g(A_t)\exp(\mu B_t)]dt$$

$$= \frac{2\lambda}{\nu^2} E[f(\exp(B_{S_\nu}))g(A_{S_\nu})\exp(\mu B_{S_\nu})],$$

where $S_\nu (= T_{\nu^2/2}$ with our previous notation) is an exponential random variable with parameter $\nu^2/2$ independent of $B$.

We apply (4.16) in the case $\mu = 0$, which we have just shown, to obtain

$$\Phi_\lambda^{(\mu)}(f,g) = \lambda \int_0^\infty \int_0^\infty g(u)f(y)y^{-(2+\nu-\mu)}p^{(\nu)}(u,1,y)dudy,$$

which implies formula (4.16) in full generality. □

The following result was used in the previous subsection.

**Theorem 4.12.** *For $\mu \in \mathbf{R}$ and $\lambda > 0$, set $\nu = \sqrt{2\lambda + \mu^2}$, $a = (\nu + \mu)/2$, $b = (\nu - \mu)/2$ and let $Z_{1,a}$ and $\gamma_b$ be a beta variable with parameters $(1, a)$ and a gamma variable with parameter $b$, respectively;*

$$P(Z_{1,a} \in dt) = a(1-t)^{a-1}dt, \qquad 0 < t < 1,$$

$$P(\gamma_b \in dt) = \frac{1}{\Gamma(b)}t^{b-1}e^{-t}dt, \qquad t > 0.$$

*Then one has the identity in law*

$$A_{T_\lambda}^{(\mu)} \stackrel{(\text{law})}{=} \frac{Z_{1,a}}{2\gamma_b}, \tag{4.18}$$

*where, on the right hand side, $Z_{1,a}$ and $\gamma_b$ are assumed to be independent.*

*Proof.* Recall the Lamperti relation (2.1) and use the same notation. Then we obtain

$$P(A^{(\mu)}_{T_\lambda} \geqq u) = \lambda \int_0^\infty e^{-\lambda t} P(A^{(\mu)}_t \geqq u) dt$$

$$= \lambda \int_0^\infty e^{-\lambda t} P(\alpha^{(\mu)}_u \leqq t) dt = E[\exp(-\lambda \alpha^{(\mu)}_u)],$$

where

$$\alpha^{(\mu)}_u = \int_0^u \frac{1}{(R^{(\mu)}_v)^2} dv.$$

Moreover, by the absolute continuity relationship (2.2) or (2.3) for the laws of different dimensional Bessel processes, we obtain

$$E[\exp(-\lambda \alpha^{(\mu)}_u)] = E[(R^{(0)}_u)^\mu \exp(-(\lambda + \mu^2/2)\alpha^{(0)}_u)] = E[(R^{(\nu)}_u)^{-2b}],$$

where $R^{(0)}_0 = R^{(\nu)}_0 = 1$.

Then, recalling formula (2.5) and using the elementary formula

$$\int_0^\infty e^{-\lambda \xi^2} \lambda^{b-1} d\lambda = \Gamma(b) \xi^{-2b},$$

we obtain after a simple change of variables

$$P(A^{(\mu)}_{T_\lambda} \geqq u) = \frac{1}{\Gamma(b)} \int_0^\infty E[\exp(-\lambda (R^{(\nu)}_u)^2)] \lambda^{b-1} d\lambda$$

$$= \frac{1}{\Gamma(b)} \int_0^{1/2u} v^{b-1} e^{-v} dv \int_{2uv}^1 a(1-t)^{a-1} dt,$$

which shows (4.18). □

## 5. Moments

In this section we present some results on the moments of $A^{(\mu)}_t$ and their consequences. We show the results on the positive moments in detail and we only consider negative moments in some special cases. We note that the expectation $E[\log(A^{(\mu)}_t)]$ has been also considered in the context of physics in Comtet-Monthus-Yor [15].

**5.1.** In this subsection we present several results about positive moments of $A_t$. At first we present an easy consequence of the Bougerol identity $\sinh(B_t) \stackrel{\text{(law)}}{=} \beta_{A_t}$.

**Proposition 5.1.** *For any $t > 0$, one has*

$$E[(A_t)^n] = \frac{1}{E[(B_1)^{2n}]} \int_{\mathbf{R}} (\sinh(x))^{2n} \frac{1}{\sqrt{2\pi t}} e^{-x^2/2t} dx, \quad n = 1, 2, ...$$

From the proposition, we obtain

$$E[(A_t)^n] = \frac{\sqrt{\pi}}{\Gamma(n+1/2)2^{3n-1}} e^{2n^2 t}(1+o(1)) \quad \text{as} \quad t \to \infty.$$

To compute $E[(A_t^{(\mu)})^n]$ for a general value of $\mu$, we note

$$E[(A_t^{(\mu)})^n] = e^{-\mu^2 t/2} E[\exp(\mu B_t)(A_t)^n],$$

which is easily shown by the Cameron-Martin theorem. For the expectation on the right hand side, we have the following (cf. [35],[65],[69]).

**Theorem 5.2.** (i) *Set $\phi(z) = z^2/2$. Then one has*

$$\int_0^\infty e^{-\lambda t} E[\exp(\mu B_t)(A_t)^n] dt = n! \prod_{j=0}^n \frac{1}{\lambda - \phi(\mu+2j)} \tag{5.1}$$

*for any $n \in \mathbf{N}$ and $\mu \geqq 0$ if $\mathrm{Re}(\lambda) > \phi(\mu+2n)$.*
(ii) *For any $t > 0, n \in \mathbf{N}, \mu \geqq 0$, it holds that*

$$E[\exp(\mu B_t)(A_t)^n] = \sum_{j=0}^n C(n,j;\mu) \exp(t\phi(\mu+2j)), \tag{5.2}$$

*where $C(0,0;\mu) = 1$ and for $j = 1,...,n$*

$$C(n,j;\mu) = 2^{-n}(-1)^{n-j} \binom{n}{j} \prod_{\substack{k \neq j \\ 0 \leqq k \leqq n}} (\mu+j+k)^{-1}.$$

*Proof.* (i) By the independent increments property of Brownian motion, we deduce

$$E[\exp(\mu B_t)(A_t)^n]$$
$$= n! E\left[\exp(\mu B_t) \int \cdots \int_{0 \leqq s_1 \cdots \leqq s_n \leqq t} \exp(2(B_{s_1}+\cdots+B_{s_n})) ds_1 \cdots ds_n\right]$$
$$= n! \int \cdots \int_{0 \leqq s_1 \cdots \leqq s_n \leqq t} E[\exp(\mu(B_t - B_{s_n}) + (\mu+2)(B_{s_n} - B_{s_{n-1}})$$
$$\qquad\qquad\qquad\qquad\qquad + \cdots + (\mu+2n)B_{s_1})] ds_1 \cdots ds_n$$
$$= n! \int \cdots \int_{0 \leqq s_1 \cdots \leqq s_n \leqq t} \exp(\phi(\mu)(t-s_n) + \phi(\mu+2)(s_n - s_{n-1})$$
$$\qquad\qquad\qquad\qquad\qquad + \cdots + \phi(\mu+2n)s_1) ds_1 \cdots ds_n.$$

Hence, setting $s_0 = 0$ and $s_{n+1} = t$, we obtain

$$\int_0^\infty e^{-\lambda t} E[\exp(\mu B_t)(A_t)^n] dt$$

$$= n! \int \cdots \int_{0 \leqq s_1 \cdots \leqq s_n \leqq s_{n+1}} \prod_{j=0}^n e^{-(\lambda - \phi(\mu+2j))(s_{j+1}-s_j)} ds_1 \cdots ds_n ds_{n+1}$$

$$= n! \prod_{j=0}^n (\lambda - \phi(\mu + 2j))^{-1}.$$

(ii) By the additive decomposition formula, we have

$$\prod_{j=0}^n \frac{1}{\lambda - \phi(\mu + 2j)} = \sum_{j=0}^n \frac{C_j^{(\mu)}}{\lambda - \phi(\mu + 2j)},$$

where

$$C_j^{(\mu)} = \prod_{\substack{k \neq j \\ 0 \leqq k \leqq n}} \frac{1}{\phi(\mu + 2j) - \phi(\mu + 2k)}.$$

Since

$$\frac{1}{\lambda - \phi(\mu + 2j)} = \int_0^\infty e^{-\lambda t} e^{\phi(\mu + 2j)t} dt,$$

we can invert the Laplace transform (5.1) and, after some elementary computations, we obtain (5.2). □

The same arguments show that formulae (5.1) and (5.2) hold when $\mu$ is any complex number. In particular, for $\mu = \sqrt{-1}\alpha$, we can show

$$E[\exp(\sqrt{-1}\alpha B_t)(A_t)^n] = \sum_{j=0}^n C(n, j; \sqrt{-1}\alpha) e^{\phi(\sqrt{-1}\alpha + 2j)t}, \quad \alpha \in \mathbf{R}. \quad (5.3)$$

By Fourier analysis, we obtain the following result on the conditional moments from (5.3). For details, see [35].

**Proposition 5.3.** *For every $t > 0$ and $n \in \mathbf{N}$, one has*

$$\frac{1}{\sqrt{2\pi t}} e^{-x^2/2t} E[(A_t)^n | B_t = x]$$
$$= \frac{e^{nx}}{n!(2\pi t^3)^{1/2}} \int_{|x|}^\infty b e^{-b^2/2t} (\cosh(b) - \cosh(|x|))^n db. \quad (5.4)$$

Recall that the Bessel function $J_0$ is defined by

$$J_0(z) = \sum_{n=0}^\infty \frac{(-1)^n (z/2)^{2n}}{(n!)^2}.$$

Then, at least heuristically, with the help of this series expansion and interchanging the order of the integration and the infinite sum on the right hand side of (5.5), the next formula holds:

$$\frac{1}{\sqrt{2\pi t}}e^{-x^2/2t}E\left[\exp\left(-\frac{1}{2}\lambda^2 A_t\right)\bigg|B_t = x\right] \\ = \frac{1}{\sqrt{2\pi t^3}}\int_{|x|}^{\infty}be^{-b^2/2t}J_0(\sqrt{2}e^{x/2}|\lambda|\sqrt{\cosh(b)-\cosh(x)})db. \quad (5.5)$$

In fact, this identity has been shown in [3] and [35] by different methods.

Comparing (5.5) with (4.4), we obtain the following expression for the heat kernel $q_\lambda(t,x,y)$ of the semi-group generated by the Schrödinger operator $H_\lambda$, $\lambda \in \mathbf{R}$, with the Liouville potential given by (4.2).

**Proposition 5.4.** *For $t > 0, x, y \in \mathbf{R}$, it holds that*

$$q_\lambda(t,x,y) = \frac{1}{\sqrt{2\pi t^3}}\int_{|x|}^{\infty}be^{-b^2/2t}J_0(\sqrt{2}e^{x/2}|\lambda|\sqrt{\cosh(b)-\cosh(x)})db. \quad (5.6)$$

*Remark* 5.1. From formula (5.6), we can show the well known explicit expression for the heat kernel of the semigroup generated by the Laplace-Beltrami operator on the hyperbolic plane. For details, see [16] and [35]. The Brownian motions on hyperbolic spaces will be discussed in the second part of the survey [51]. See also [30] and [47] about the close relationship between the Brownian motions on hyperbolic spaces and our exponential functionals.

*Remark* 5.2. Denoting by $C_n$ the right hand side of formula (5.4), it is easy to show

$$\sum_{n=1}^{\infty}(C_n)^{-1/2n} < \infty.$$

This means that Carleman's sufficient condition for the unique solvability of the Stieltjes moment problem does not hold for the (conditional) distribution of $A_t$. Recently Hörfelt [34] announced that he has proven the indeterminacy for the moment problem for $A_t$.

Dufresne [22] has shown that the Stieltjes moment problem for the reciprocal $(A_t^{(\mu)})^{-1}$ is determinate. By using this fact and formula (4.1), he has also shown another representation for the density of $A_t^{(\mu)}$ after giving another expression for the moments. See also [6] for a similar result for exponential functionals of Lévy processes.

We go back to formula (5.4). Let us consider an exponential random variable **e** whose density is $e^{-x}$, $x \geqq 0$. Then we have $E[\mathbf{e}^n] = n!$ and

$$P(\sqrt{2\mathbf{e}t} > b|\sqrt{2\mathbf{e}t} > |x|) = e^{x^2/2t}\int_b^{\infty}\frac{c}{t}e^{-c^2/2t}dc, \quad b > |x|.$$

Hence, assuming that **e** and $B$ are independent, we obtain from (5.4)

$$E[(e^{-x}\mathbf{e}A_t)^n|B_t = x] = E[(\cosh(\sqrt{2\mathbf{e}t})-\cosh(x))^n|\sqrt{2\mathbf{e}t} > |x|].$$

Although we do not know the unique solvability of the moment problem for the conditional distribution of $A_t$, we can show the following result. It says that, if we multiply $A_t$ by an independent exponential random variable, the probability law of $A_t$ is modified into a rather simple one.

**Theorem 5.5.** *For any $t > 0$, the conditional distribution of $\mathbf{e}A_t$ conditionally on $B_t = x$ coincides with that of $e^x(\cosh(\sqrt{2\mathbf{e}t}) - \cosh(x))$ given $\mathbf{e} > x^2/2t$.*

For a proof, see [20] or [48].

In these articles one also finds the following results, which are shown by using the expression (4.1) for the joint density of $(A_t, B_t)$.

**Theorem 5.6.** (i) *For any $t > 0$ and $\lambda > 0$, one has*

$$E\left[\exp\left(-\frac{\lambda}{A_t}\right)\bigg|B_t = x\right] = \exp\left(-\frac{\varphi_x(\lambda)^2 - x^2}{2t}\right), \tag{5.7}$$

*where $\varphi_x(\lambda) = \mathrm{Argcosh}(\lambda e^{-x} + \cosh(x))$.*
(ii) *Letting $\mathbf{e}$ be an exponential random variable independent of $B$ and $\{L_u, u \geqq 0\}$ be the local time of $B$ at $0$, one has the identity in law for fixed $t$*

$$(\mathbf{e}A_t, B_t) \stackrel{(\mathrm{law})}{=} (e^{B_t}(\cosh(|B_t| + L_t) - \cosh(B_t)), B_t).$$

*Proof.* We only give a proof of formula (5.7), which will be used in the next subsection. Let $f$ be a non-negative Borel function on $\mathbf{R}$. Then, from formula (4.1), we deduce

$$\int_{\mathbf{R}} f(x) E\left[\exp\left(-\frac{\lambda}{A_t}\right)\bigg|B_t = x\right]\frac{1}{\sqrt{2\pi t}}e^{-x^2/2t}dx$$
$$= E\left[f(B_t)\exp\left(-\frac{\lambda}{A_t}\right)\right]$$
$$= \int_{\mathbf{R}} f(x)dx \int_0^\infty \exp\left(-\frac{\lambda + e^x\cosh(x)}{u}\right)\theta(e^x/u, t)\frac{du}{u}$$
$$= \int_{\mathbf{R}} f(x)dx \int_0^\infty \exp(-\varphi_x(\lambda)w)\theta(w, t)\frac{dw}{w}.$$

Then, with the help of formula (4.6), we obtain

$$\int_{\mathbf{R}} f(x) E\left[\exp\left(-\frac{\lambda}{A_t}\right)\bigg|B_t = x\right]\frac{1}{\sqrt{2\pi t}}e^{-x^2/2t}dx$$
$$= \int_{\mathbf{R}} f(x)\frac{1}{\sqrt{2\pi t}}\exp\left(-\frac{\varphi_x(\lambda)^2}{2t}\right)dx.$$

$\square$

**5.2.** As is mentioned in Remark 5.2, it is important to have formulae for the negative moments for $A_t^{(\mu)}$. Moreover, considering the negative (conditional) moments of $A_t^{(\mu)}$, we obtain several interesting identities which are related to some special functions. However, we restrict ourselves to some special cases which can be obtained from the results mentioned in the previous sections. For more results, see [3], [19], [20], [22] and the references cited therein.

We again start from the Bougerol identity or rather from formula (3.1).

**Proposition 5.7.** *For any $t > 0$ and $\alpha \in \mathbf{R}$, one has*

$$E[(A_t)^{-1/2}] = t^{-1/2} \quad and \quad E\left[\left(\int_0^1 \exp(\alpha B_u)du\right)^{-1/2}\right] = 1.$$

The first identity is a direct consequence of formula (3.1), and, from this, we obtain the second striking one by the scaling property of Brownian motion.

The following is also a consequence of the Bougerol identity.

**Proposition 5.8.** *For any $p$ with $0 < p < 1/2$, one has*

$$E[(A_t)^{-p}] = \frac{\Gamma(\alpha)}{2^{1/2-p}\sqrt{\pi}} t^{-1/2}(1 + o(1)) \quad as \quad t \to \infty. \tag{5.8}$$

In fact, Hariya-Yor [32] have shown that formula (5.8) holds for any $p > 0$. In [32] one finds a complete description on the asymptotic behavior of $E[(A_t^{(\mu)})^p]$, $p \in \mathbf{R}$, as $t \to \infty$ and on the limit of the probability measure

$$(E[(A_t^{(\mu)})^p])^{-1}(A_t^{(\mu)})^p dP.$$

We will present a part of the results of [32] in Section 7.

To find an analogous result to that of Proposition 5.7 for the conditional distribution of $A_t$ given $B_t$, we recall formula (5.7). By differentiating both hand sides of (5.7) in $\lambda$, we can compute $E[(A_t)^{-n}|B_t = x]$, $n = 1, 2, ...$ In particular, we obtain the following.

**Proposition 5.9.** *For any $t > 0$ and $x \in \mathbf{R}$, one has*

$$E[(A_t)^{-1}|B_t = x] = \frac{xe^{-x}}{t\sinh(x)} \quad and \quad E[(A_t)^{-1}|B_t = 0] = t^{-1}.$$

By the scaling property of Brownian motion, we also obtain

$$E\left[\left(\int_0^1 \exp(\alpha B_u)du\right)^{-1}\bigg|B_1 = 0\right] = 1 \quad \text{for any } \alpha \in \mathbf{R}.$$

Chaumont-Hobson-Yor [13] have shown that this striking identity is a consequence of the so-called cyclic exchangeability property of the Brownian bridge and they have obtained several similar identities for the stochastic processes which enjoy this property.

## 6. Perpetual functionals

Throughout this section we assume $\mu > 0$. For a non-negative Borel measurable function $f$, we set

$$A_\infty^{(\mu)}(f) = \int_0^\infty f(B_s^{(\mu)})ds$$

and call it a perpetual functional.

At first we show a necessary and sufficient condition in order that $A_\infty^{(\mu)}(f)$ is finite almost surely. We follow Salminen-Yor [59]. See also Engelbert-Senf [24] for another approach.

**Theorem 6.1.** *Let $\mu$ be positive and $f$ be a non-negative and locally integrable Borel function on $\mathbf{R}$. Then, in order for $A_\infty^{(\mu)}(f)$ to be finite almost surely, it is necessary and sufficient that $\int^\infty f(y)dy < \infty$.*

Before proceeding to a proof, we recall a lemma due to Jeulin [38].

**Lemma 6.1.** *Let $\{Z_y, y \geqq 0\}$ be a stochastic process such that the probability law of $Z_y$ does not depend on $y$ and is absolutely continuous with respect to the Lebesgue measure. Assume also that $E[|Z_y|] < \infty$. Then, for any $\sigma$-finite non-negative measure $\nu$ on $\mathbf{R}_+$, $\int_0^\infty |Z_y|\nu(dy)$ is finite almost surely if and only if $\nu(\mathbf{R}_+) < \infty$.*

*Proof of Theorem 6.1.* Letting $L^{(\mu)}(t,y)$ be the local time of $B^{(\mu)}$ at $y$ up to time $t$, we have

$$\int_0^t f(B_s^{(\mu)})ds = \int_\mathbf{R} f(y)L^{(\mu)}(t,y)dy$$

and, letting $t$ tend to $\infty$ and using the monotone convergence theorem, we obtain

$$A_\infty^{(\mu)}(f) = \int_\mathbf{R} f(y)L^{(\mu)}(\infty, y)dy.$$

The distribution of $L^{(\mu)}(\infty, y)$ is known (cf. [59]) and is given by

$$P(L^{(\mu)}(\infty, y) \geqq v) = \begin{cases} e^{-\mu v}, & y \geqq 0, \\ e^{2\mu y - \mu v}, & y < 0. \end{cases}$$

The probability law of $\{L^{(\mu)}(\infty, y), y \geqq 0\}$ does not depend on $y$. Hence, letting $\widehat{f}$ be the non-negative function given by $\widehat{f}(x) = f(x)\mathbf{1}_{\mathbf{R}_+}(x)$, we obtain from Lemma 6.1

$$P\left(\int_0^\infty \widehat{f}(B_s^{(\mu)})ds < \infty\right) = 1 \quad \text{if and only if} \quad \int^\infty \widehat{f}(x)dx < \infty.$$

The rest of the proof is easy. □

In [59], sufficient conditions on $f$ for $A_\infty^{(\mu)}(f)$ to have moments of any order and to have exponential moments are also given.

We next consider the important and typical case, $f(x) = e^{-2x}$. We set

$$A_\infty^{(-\mu)} = \int_0^\infty \exp(2B_s^{(-\mu)})ds \left( \stackrel{(\text{law})}{=} \int_0^\infty \exp(-2B_s^{(\mu)})ds \right).$$

The following Theorem 6.2 is due to Dufresne [21].

Bertoin-Yor [6] (see also [7]) have considered Lévy processes which drift to $\infty$ as $t \to \infty$ and have shown that the distributions of the corresponding perpetual functionals are determined by their negative moments. Moreover, they have given another proof of the following theorem. See also [57], p.452, and [70].

**Theorem 6.2.** *Let $\mu > 0$. Then, $A_\infty^{(-\mu)}$ is distributed as $(2\gamma_\mu)^{-1}$, where $\gamma_\mu$ denotes a gamma random variable with parameter $\mu$.*

*Proof.* Differently from the proofs given in the above cited references and at the end of this section, we give an analytic proof.

Let $\{e_x^{(\mu)}(t)\}$ be the diffusion process given by $e_x^{(\mu)}(t) = x\exp(B_t^{(\mu)})$ and $\tau_z$ be its first hitting time of $z$. For $\alpha > 0$, we set

$$v_z(x; \alpha) = E\left[\exp\left(-\frac{\alpha^2}{2}\int_0^{\tau_z}(e_x^{(\mu)}(t))^{-2}dt\right)\right].$$

Then $v_z$ solves the equation

$$\frac{1}{2}x^2 v''(x) + \left(\mu + \frac{1}{2}\right)xv'(x) = \frac{\alpha^2}{2x^2}v(x),$$
$$v(z) = 1, \qquad \lim_{x\downarrow 0}v(x) = 0.$$

Moreover, we have

$$\begin{aligned}\lim_{z\to\infty} v_z(x;\alpha) &= E\left[\exp\left(-\frac{\alpha^2}{2x^2}\int_0^\infty e^{-2B_t^{(\mu)}}dt\right)\right] \\ &= E\left[\exp\left(-\frac{\alpha^2}{2x^2}A_\infty^{(-\mu)}\right)\right].\end{aligned} \qquad (6.1)$$

Since $v_z(x;\alpha)$ depends only on $x/\alpha$, we denote it by $\widetilde{v}_z(x/\alpha)$.

On the other hand, making the change of variables $\xi = \alpha/x$ and $\widetilde{v}_z(x/\alpha) = \xi^\mu \phi(\xi)$, we obtain after an elementary computation

$$\phi''(\xi) + \frac{1}{\xi}\phi'(\xi) - \left(1 + \frac{\mu^2}{\xi^2}\right)\phi(\xi) = 0.$$

This is the differential equation for the modified Bessel functions $I_\mu$ and $K_\mu$ of order $\mu$. Since $\phi(\xi) \to 0$ as $\xi \to \infty$ by the boundary condition for $v_z$, we easily obtain

$$\phi(\xi) = \left(\frac{z}{\alpha}\right)^\mu \frac{K_\mu(\xi)}{K_\mu(\alpha/z)}$$

and, therefore,
$$\widetilde{v}_z(x/\alpha) = \left(\frac{z}{x}\right)^\mu \frac{K_\mu(\alpha/x)}{K_\mu(\alpha/z)}.$$

Using the formula
$$K_\mu(\eta) = 2^{\mu-1}\Gamma(\mu)\eta^{-\mu}(1+o(1)) \quad \text{as} \quad \eta \to 0$$

and setting $x = 1$, we get from (6.1)
$$E\left[\exp\left(-\frac{1}{2}\alpha^2 A_\infty^{(-\mu)}\right)\right] = \frac{2}{\Gamma(\mu)}\left(\frac{\alpha}{2}\right)^\mu K_\mu(\alpha). \tag{6.2}$$

Finally, with the help of the integral representation for $K_\mu$ (cf. [41], p.119)
$$K_\mu(\alpha) = \frac{1}{2}\left(\frac{\alpha}{2}\right)^\mu \int_0^\infty e^{-\xi-\alpha^2/4\xi}\xi^{-\mu-1}d\xi$$

and the injectivity of the Laplace transform, we obtain the assertion. $\square$

Theorem 6.2 has been extended in several ways. For example, in [58], the following has been shown.

**Proposition 6.3.** *Setting*
$$A_\infty^{(-\mu,\pm)} = \int_0^\infty e^{2B_s^{(-\mu)}} \mathbf{1}_{\{B_s^{(-\mu)} \in \mathbf{R}_\pm\}} ds,$$

*one has*
$$E\left[\exp\left(-\frac{1}{2}(\alpha^2 A_\infty^{(-\mu,+)} + \beta^2 A_\infty^{(-\mu,-)})\right)\right]$$
$$= \frac{\beta^\mu}{\Gamma(\mu)2^{\mu-1}} \frac{K_\mu(\alpha)}{\beta K_\mu(\alpha)I_{\mu-1}(\beta) + \alpha K_{\mu-1}(\alpha)I_\mu(\beta)} \tag{6.3}$$

*for any $\alpha, \beta > 0$.*

We recover formula (6.2) from (6.3) by taking $\alpha = \beta$ and using the formula
$$K_\mu(z)I_{\mu-1}(z) + K_{\mu-1}(z)I_\mu(z) = \frac{1}{z}.$$

Moreover, by the Lamperti relation, we have
$$A_\infty^{(-\mu)} = \inf\{u; R_u^{(-\mu)} = 0\},$$

where $\{R_u^{(-\mu)}\}$ is a Bessel process with index $-\mu < 0$ starting from 1. Hence, by standard diffusion theory, we can verify
$$E\left[\exp\left(-\frac{1}{2}\alpha^2 A_\infty^{(-\mu)}\right)\right] = \frac{2}{\Gamma(\mu)}\left(\frac{\alpha}{2}\right)^\mu K_\mu(\alpha)$$

and we have another proof for Theorem 6.2.

In Salminen-Yor [60], the distributions of several perpetual functionals have been identified with those of the hitting and occupation times of Brownian motion and Bessel processes. See also [10] for more results.

When $\mu \geqq 0$, $A_t^{(\mu)} \to \infty$ as $t \to \infty$. We now show the following limit theorem for $\log(A_t^{(\mu)})$.

**Proposition 6.4.** *Let $Z$ be a standard normal random variable. Then, as $t \to \infty$ :*
*(i) if $\mu = 0$, $t^{-1/2} \log(A_t^{(0)})$ converges in law to $2|Z|$;*
*(ii) if $\mu > 0$, $t^{-1/2}(\log(A_t^{(\mu)}) - 2\mu t)$ converges in law to $2Z$.*

*Proof.* (i) By the scaling property of Brownian motion, we have

$$A_t^{(0)} \stackrel{\text{(law)}}{=} t \int_0^1 e^{2\sqrt{t} B_s} ds$$

and, by the Laplace method,

$$\frac{1}{\sqrt{t}} \log\left(\int_0^1 e^{2\sqrt{t} B_s} ds\right) \to 2 \max_{0 \leqq s \leqq 1} B_s.$$

Since $\max_{0 \leqq s \leqq 1} B_s \stackrel{\text{(law)}}{=} |Z|$, we obtain the result. Note that another proof may easily be obtained with the help of Bougerol's identity.
(ii) By the invariance of the probability law of Brownian motion by time reversal from a fixed time, we have

$$A_t^{(\mu)} \stackrel{\text{(law)}}{=} e^{2B_t^{(\mu)}} \int_0^t e^{-2B_s^{(\mu)}} ds$$

hence

$$\frac{1}{\sqrt{t}}\left(\log(A_t^{(\mu)}) - 2\mu t\right) \stackrel{\text{(law)}}{=} 2\frac{B_t}{\sqrt{t}} + \frac{1}{\sqrt{t}} \log\left(\int_0^t e^{-2B_s^{(\mu)}} ds\right).$$

Since $t^{-1/2} B_t$ is identical in law with $Z$ and the second term on the right hand side tends to 0 as $t \to \infty$ almost surely, we obtain the result. □

By using the Markov property of Brownian motion, we can show the following process level extension of Proposition 6.4.

**Proposition 6.5.** *In the sense of convergence of finite dimensional distributions,*
*(i) if $\mu = 0$, $\{\lambda^{-1/2} \log(A_{\lambda t}^{(0)}), t > 0\}$ converges to $\{2\max_{0 \leqq s \leqq t} B_s, t > 0\}$ as $\lambda \to \infty$ ;*
*(ii) if $\mu > 0$, $\{\lambda^{-1/2}(\log(A_{\lambda t}^{(\mu)}) - 2\mu\lambda t), t > 0\}$ converges to $\{2B_t, t > 0\}$.*

We end this section by presenting some results on the asymptotics of the joint distribution of $B_t$ and $A_t$ as $t \to \infty$.

In Section 2 we showed

$$\lim_{t\to\infty} \sqrt{2\pi t^3}\theta(r,t) = K_0(r), \quad r>0.$$

From this result, we can show the following.

**Proposition 6.6.** *For any non-negative functions $f$ on $\mathbf{R}$ and $g$ on $\mathbf{R}_+$, satisfying*

$$\int_{-\infty}^{\infty} |x|f(x)dx < \infty, \quad \int^{\infty} f(x)dx < \infty, \quad \int_0^{\infty} g(y)e^{-1/2y}\frac{dy}{y} < \infty,$$

*one has*

$$\lim_{t\to\infty} \sqrt{2\pi t^3} E[f(B_t)g(A_t)] = \int_{\mathbf{R}}\int_0^{\infty} f(x)g(y)K_0(e^x/y)\exp\left(-\frac{1+e^{2x}}{2y}\right)\frac{dxdy}{y}.$$

This result has been extended by Kotani [39] and, recently, by Hariya [31] to general Brownian functionals of the form $\int_0^t V(B_s)ds$ for a non-negative function $V$, instead of the exponential function. They have shown, under suitable assumptions, that

$$\lim_{t\to\infty} t^{3/2} E\left[f(B_t)g\left(\int_0^t V(B_s)ds\right)\right]$$

exists and may be identified.

Kotani used Krein's theory of strings and Hariya obtained the result by reducing the computations to those for the three-dimensional Bessel process. For details, see the original papers.

*Remark* 6.1. As is mentioned above, the exponential perpetual functionals of Lévy processes have been studied in [6] and [7].

On the other hand, for two independent Lévy processes $\{Y_t\}$ and $\{Z_t\}$, which are given by the sums of constant drifts, Brownian motions and compound Poisson processes, the perpetual functional $\int_0^{\infty} \exp(Y_t)dZ_t$ has been studied in Gjessing-Paulsen [27], Nilsen-Paulsen [54], Paulsen [56] and so on. For details, see the original papers.

## 7. Some limiting distributions

In this section, following Hariya-Yor [32], we present the results on the asymptotic behavior of

$$D_t^{(\mu,\alpha)} = E[\exp(-\alpha A_t^{(\mu)})] \quad \text{and} \quad \Delta_t^{(\mu,m)} = E[(A_t^{(\mu)})^{-m}]$$

as $t\to\infty$ and on the weak limits of the probability measures

$$(D_t^{(\mu,\alpha)})^{-1}\exp(-\alpha A_t^{(\mu)})dP \quad \text{and} \quad (\Delta_t^{(\mu,m)})^{-1}(A_t^{(\mu)})^{-m}dP.$$

We present the results without proofs. For details and proofs, see the original paper.

**Theorem 7.1.** *Fix $\alpha > 0$. Then we have the following:*
(i) *If $\mu > 0$, it holds that*

$$\lim_{t\to\infty} t^{3/2} e^{\mu^2 t/2} D_t^{(\mu,\alpha)} = \frac{2^{\mu-3/2}}{\sqrt{\pi}} \left(\Gamma\left(\frac{\mu}{2}\right)\right)^2 \frac{K_0(\sqrt{2\alpha})}{(\sqrt{2\alpha})^\mu}.$$

(ii) *If $\mu = 0$, it holds that*

$$\lim_{t\to\infty} t^{1/2} D_t^{(0,\alpha)} = \sqrt{\frac{2}{\pi}} K_0(\sqrt{2\alpha}).$$

(iii) *If $\mu < 0$, it holds that*

$$\lim_{t\to\infty} D_t^{(\mu,\alpha)} = \frac{2^{\mu+1}}{\Gamma(|\mu|)} \frac{K_\mu(\sqrt{2\alpha})}{(\sqrt{2\alpha})^\mu}.$$

For $\Delta_t^{(\mu,m)}$, we consider the expectation of a more general functional. For $a \geqq 0, \xi > 0$, we set

$$\Delta_t^{(\mu,m)}(a,\xi) = E[(a + \xi A_t^{(\mu)})^{-m}].$$

To present the result, we need to partition the $(\mu, m)$-plane into the following six subsets:

$$\begin{aligned}
R_1 &= \{2m > \mu, \mu > 0\}, & L_1 &= \{2m = \mu, \mu > 0\}, \\
R_2 &= \{m < \mu, 2m < \mu\}, & L_2 &= \{m = \mu, \mu < 0\}, \\
R_3 &= \{m > \mu, \mu < 0\}, & L_3 &= \{m > 0, \mu = 0\}.
\end{aligned}$$

Here are some corresponding limiting results.

**Theorem 7.2.** (i) *If $(\mu, m) \in R_1$, one has*

$$\lim_{t\to\infty} t^{3/2} e^{\mu^2 t/2} \Delta_t^{(\mu,m)}(a,\xi)$$
$$= \frac{2^{(\mu-5)/2}}{\sqrt{\pi}} \Gamma\left(\frac{\mu}{2}\right) B\left(\frac{\mu}{2}, m - \frac{\mu}{2}\right) \xi^{-\mu/2} \int_0^\infty \frac{x^{m-1-\mu/2} e^{-x/2}}{(ax+\xi)^{m-\mu/2}} dx.$$

(ii) *If $(\mu, m) \in L_1$, one has*

$$\lim_{t\to\infty} t^{1/2} e^{\mu^2 t/2} \Delta_t^{(\mu,m)}(a,\xi) = \frac{2^{m-1/2}}{\sqrt{\pi}} \Gamma(m) \xi^{-m}.$$

(iii) *If $(\mu, m) \in R_2$, one has*

$$\lim_{t\to\infty} e^{2m(\mu-m)t} \Delta_t^{(\mu,m)}(a,\xi) = \frac{2^m \Gamma(\mu-m)}{\Gamma(\mu-2m)} \xi^{-m}.$$

(iv) If $(\mu, m) \in L_2$, one has

$$\lim_{t \to \infty} t^{-1} \Delta_t^{(\mu,m)}(a, \xi) = \frac{2^{m+1}|\mu|}{\Gamma(|\mu|)} \xi^{-m}.$$

(v) If $(\mu, m) \in R_3$, one has

$$\lim_{t \to \infty} \Delta_t^{(\mu,m)}(a, \xi) = \frac{2^\mu}{\Gamma(|\mu|)} \int_0^\infty \frac{x^{m-\mu-1} e^{-x/2}}{(ax+\xi)^m} dx.$$

(vi) If $(\mu, m) \in L_3$, one has

$$\lim_{t \to \infty} t^{1/2} \Delta_t^{(\mu,m)}(a, \xi) = \frac{1}{\sqrt{2\pi}} \int_0^\infty \frac{x^{m-1} e^{-x/2}}{(ax+\xi)^m} dx.$$

Next, letting $\mathcal{W}^{(\mu)}$ be the probability law of $B^{(\mu)}$ on the canonical path space $W = C([0,\infty) \to \mathbf{R})$, we consider the Gibbsian-like measure $\mathbb{P}_t^{(\mu,\alpha)}$ and the moment density measure $\mathcal{P}_t^{(\mu,m)}$ on $W$ defined by

$$d\mathbb{P}_t^{(\mu,\alpha)}(X) = (D_t^{(\mu,\alpha)})^{-1} \exp(-\alpha A_t(X)) d\mathcal{W}^{(\mu)}(X)$$

and

$$d\mathcal{P}_t^{(\mu,m)}(X) = (\Delta_t^{(\mu,m)})^{-1} (A_t(X))^{-m} d\mathcal{W}^{(\mu)}(X)$$

for $\mu \in \mathbf{R}, \alpha > 0$ and $m \in \mathbf{R}$, where $A_t(X) = \int_0^t \exp(2X_s) ds$.

One of the aims in [32] is to show that the probability measures $\mathbb{P}_t^{(\mu,\alpha)}$ and $\mathcal{P}_t^{(\mu,m)}$ converge weakly on $W$ as $t \to \infty$ and to characterize the limiting measures.

To present the results, we need to recall the generalized inverse Gaussian distribution, which we denote by $\mathrm{GIG}(\mu; a, b)$, and its relation to Brownian motion with drift.

For $a, b > 0$ and $\mu \in \mathbf{R}$, the generalized inverse Gaussian distribution is the probability distribution on $\mathbf{R}_+$ whose density is given by

$$\left(\frac{b}{a}\right)^\mu \frac{x^{\mu-1}}{2K_\mu(ab)} \exp\left(-\frac{1}{2}\left(b^2 x + \frac{a^2}{x}\right)\right) dx, \quad x > 0.$$

In particular, $\mathrm{GIG}(-1/2; a, b)$ is called the inverse Gaussian distribution.

Let $\mu > 0$ and denote by $\tau_{x \to y}^{(\mu)}$ and $\ell_{x \to y}^{(\mu)}$ the first hitting time at $y$ and the last exit time from $y$ of $\{x + B^{(\mu)}(t)\}$, respectively. Then, it is well known that the distribution of $\tau_{x \to y}^{(\mu)}$ is $\mathrm{GIG}(-1/2; y-x, \mu)$ and that of $\ell_{x \to y}^{(\mu)}$ is $\mathrm{GIG}(1/2; y-x, \mu)$. When $x = y$, $\ell_{y \to y}^{(\mu)}$ is distributed as $2\mu^{-2} \gamma_{1/2}$ for a gamma random variable $\gamma_{1/2}$ with parameter $1/2$. By the strong Markov property of Brownian motion, we have for $x < y$

$$\ell_{x \to y}^{(\mu)} \overset{(\text{law})}{=} \tau_{x \to y}^{(\mu)} + \ell_{y \to y}^{(\mu)},$$

which yields the special case of the following identity in law (7.1) (with $\nu = 1/2$, $a = y - x$, $b = \mu$).

Generally, we let $I_{a,b}^{(\pm\nu)}$ be GIG$(\pm\nu; a, b)$ random variables and $\gamma_\nu$ be a gamma variable with parameter $\nu$. Then we have

$$I_{a,b}^{(\nu)} \stackrel{(\text{law})}{=} I_{a,b}^{(-\nu)} + \frac{2}{b^2}\gamma_\nu \qquad (7.1)$$

for all $\nu > 0$, where we assume that the two random variables on the right hand side are independent. We can easily prove this identity by computing the Laplace transforms. For more details of this identity, see [43] and [62].

For our purpose, we only need to discuss a particular case and let $I(\mu, a)$ be a GIG$(\mu; a, 1)$ random variable.

Moreover, letting $T_z, z > 0$, be the transformation on the path space $W$ given by

$$T_z(X)(t) = X_t - \log(1 + zA_t(X))$$

and assuming that $I(\mu, a)$ is independent of $B^{(\mu)}$, we denote by $\mathcal{W}_{I(\mu,a)}^{(\mu)}$ the probability law of the stochastic process $\{T_{I(\mu,a)}(B^{(\mu)})(t), t \geqq 0\}$:

$$\mathcal{W}_{I(\mu,a)}^{(\mu)} = P \circ (T_{I(\mu,a)}(B^{(\mu)}))^{-1}.$$

**Theorem 7.3.** (i) *If $\mu \geqq 0$, then $\mathbb{P}_t^{(\mu,\alpha)}$ converges weakly to $\mathcal{W}_{I(0,\sqrt{2\alpha})}^{(0)}$ on $W$ as $t \to \infty$.*
(ii) *If $\mu < 0$, then $\mathbb{P}_t^{(\mu,\alpha)}$ converges weakly to $\mathcal{W}_{I(-\mu,\sqrt{2\alpha})}^{(-\mu)}$, which is equal to $\mathcal{W}_{I(\mu,\sqrt{2\alpha})}^{(\mu)}$.*

Note that, if $\mu \geqq 0$, the limiting measure does not depend on $\mu$. Some heuristic discussions and analytical explanations are found in the original paper [32].

For the moment density measures $\mathcal{P}_t^{(\mu,m)}$, we have the following results.

**Theorem 7.4.** *Let $\gamma_\nu$ be a gamma random variable with parameter $\nu > 0$ independent of $B^{(\mu)}$. Then, the following holds:*
(i) *If $(\mu, m) \in R_1 \cup L_3$, $\mathcal{P}_t^{(\mu,m)}$ converges weakly to $\mathcal{W}_{2\gamma_{m-\mu/2}}^{(0)}$ on $W$ as $t \to \infty$;*
(ii) *If $(\mu, m) \in L_1 \cup R_2 \cup L_2$, $\mathcal{P}_t^{(\mu,m)}$ converges weakly to $\mathcal{W}^{(\mu-2m)}$;*
(iii) *If $(\mu, m) \in R_3$, $\mathcal{P}_t^{(\mu,m)}$ converges weakly to $\mathcal{W}_{2\gamma_{m-\mu}}^{(-\mu)}$.*

*Remark* 7.1. Both the generalized inverse Gaussian distributions and the path transformations $T_z$ will also play important roles when we study some diffusion processes associated with the exponential functionals. For details, see [18, 50]. This will be discussed in the second part of our survey [51].

## 8. A 3*D* distribution

In the different studies of the so-called stochastic volatility models in mathematical finance or in those of generalized Laplacians on the hyperbolic plane and usual Laplacians on complex hyperbolic spaces, or again in those of some diffusion processes in random environments (see, e.g., [16], [35], [42], [45], [46]), one often needs to consider, together with $A_t^{(\mu)}$, another exponential functional

$$a_t^{(\mu)} = \int_0^t \exp(B_s^{(\mu)})ds$$

and to study the joint distribution of $(a_t^{(\mu)}, A_t^{(\mu)}, B_t^{(\mu)})$.

In this section, following [1], we show that our method to prove Theorem 4.1 may also be applied to this study and that we can give an explicit expression for the Laplace transform of the conditional distribution of $A_t^{(\mu)}$ given $(a_t^{(\mu)}, B_t^{(\mu)})$.

To present the result, we denote by $\psi_t^{(\mu)}(v,x)$ the joint density of $(a_t^{(\mu)}, B_t^{(\mu)})$ given by

$$\psi_t^{(\mu)}(v,x) = \frac{1}{16} e^{\mu x - \mu^2 t/2} v^{-1} \exp\left(-\frac{2(1+e^x)}{v}\right) \theta(4e^{x/2}/v, t/4),$$

which is easily shown from (4.1) by the scaling property of Brownian motion. We also set

$$\phi \equiv \phi(v,x;\lambda) = \frac{2\lambda \exp(x/2)}{\sinh(\lambda v/2)}.$$

**Theorem 8.1.** *For any $t > 0, \lambda > 0, v > 0$ and $x \in \mathbf{R}$, one has*

$$E\left[\exp\left(-\frac{1}{2}\lambda^2 A_t^{(\mu)}\right) \bigg| a_t^{(\mu)} = v, B_t^{(\mu)} = x\right] \psi_t^{(\mu)}(v,x)$$

$$= e^{\mu x - \mu^2 t/2} \frac{\lambda}{4\sinh(\lambda v/2)} \exp(-\lambda(1+e^x)\coth(\lambda v/2)) \theta(\phi, t/4).$$

*Proof.* We only give a proof when $\mu = 0$. Once this is done, the result for a general value of $\mu$ is easily obtained by the Cameron-Martin theorem. We denote $A_t^{(0)}$ and $a_t^{(0)}$ simply by $A_t$ and $a_t$, respectively.

For $\lambda > 0$ and $k \in \mathbf{R}$, we consider the following Schrödinger operator $H_{\lambda,k}$ with the Morse potential $V_{\lambda,k}$ on $\mathbf{R}$:

$$H_{\lambda,k} = -\frac{1}{2}\frac{d^2}{dx^2} + V_{\lambda,k}, \qquad V_{\lambda,k}(x) = \frac{1}{2}\lambda^2 e^{2x} - \lambda k e^x.$$

Letting $q_{\lambda,k}(t,x,y)$ be the heat kernel of the semigroup generated by $H_{\lambda,k}$, we have by the Feynman-Kac formula

$$q_{\lambda,k}(t,x,y)$$
$$= E\left[\exp\left(-\frac{1}{2}\lambda^2 e^{2x} A_t + \frac{1}{2}\lambda k e^x a_t\right) \bigg| B_t = y - x\right] \frac{1}{\sqrt{2\pi t}} e^{-|y-x|^2/2t}.$$

The Green function $G_{\lambda,k}$ for $H_{\lambda,k}$ is obtained by using the general theory of the Sturm-Liouville operators. In fact, letting $W_{k,\alpha}$ and $M_{k,\alpha}$ be the Whittaker functions (cf. [12], [41]), we can show

$$G_{\lambda,k}(x,y;\alpha^2/2) = \int_0^\infty e^{-\alpha^2 t/2} q_{\lambda,k}(t,x,y) dt$$
$$= \frac{\Gamma(\alpha - k + 1/2)}{\lambda \Gamma(1 + 2\alpha)} e^{-(x+y)/2} W_{k,\alpha}(2\lambda e^y) M_{k,\alpha}(2\lambda e^x),$$

for $y \geqq x$ and $\alpha > 0$, where, if $k > 0$, we assume $\alpha > k - 1/2$. See also [9] p.196.

Note the integral representation for the product of the Whittaker functions (see [12], p.86, and [29], p.729):

$$W_{k,\alpha}(\xi) M_{k,\alpha}(\eta)$$
$$= \frac{\sqrt{\xi\eta}\Gamma(1+2\alpha)}{\Gamma(\alpha-k+1/2)} \int_0^\infty e^{-(\xi+\eta)\cosh(z)/2} \left(\coth\left(\frac{z}{2}\right)\right)^{2k} I_{2\alpha}(\sqrt{\xi\eta}\sinh(z)) dz.$$

Then we obtain

$$G_{\lambda,k}(x,y;\alpha^2/2) = 2\int_0^\infty \frac{e^{2ku}}{\sinh(u)} e^{-\lambda(e^x+e^y)\coth(u)} I_{2\alpha}\left(\frac{2\lambda e^{(x+y)/2}}{\sinh(u)}\right) du$$

by a simple change of variables, $e^u = \coth(z/2)$.

Recalling the characterization (2.10) for the function $\theta(r,t)$, we obtain

$$q_{\lambda,k}(t,x,y) = \int_0^\infty \frac{e^{2ku}}{2\sinh(u)} e^{-\lambda(e^x+e^y)\coth(u)} \theta(\bar{\phi}, t/4) du$$

by the injectivity of Laplace transform, where $\bar{\phi} = 2\lambda e^{(x+y)/2}/\sinh(u)$. The rest of proof is easy and is omitted. □

Another probabilistic proof of the theorem is given in [1].

## Appendix A: Density of the Hartman-Watson distributions

In this section, following [67], we show that the density of the Hartman-Watson distribution is given by $(I_0(r))^{-1}\theta(r,t)$, where the function $\theta(r,t)$ is given by (2.9).

**Theorem A.1.** *Let $\theta(r,\cdot), r > 0$, be the function on $(0,\infty)$ given by (2.9). Then we have*

$$\int_0^\infty e^{-\nu^2 t/2} \theta(r,t) dt = I_\nu(r), \qquad \nu > 0, r > 0.$$

*Proof.* We start from the following well known integral representation for the modified Bessel function $I_\nu(z)$ (see, e.g., [64], p.176):

$$I_\nu(z) = \frac{1}{2\pi\sqrt{-1}} \int_C e^{z\cosh(w) - \nu w} dw, \qquad \text{Re}(z) > 0, \text{Re}(\nu) > -1/2,$$

where the path $C$ of integration consists of three sides of the infinite rectangle with vertices at $\infty - \pi\sqrt{-1}, -\pi\sqrt{-1}, \pi\sqrt{-1}$ and $\infty + \pi\sqrt{-1}$. Using the elementary formula

$$e^{-\nu w} = \int_0^\infty e^{-\nu^2 t/2} \frac{w}{\sqrt{2\pi t^3}} e^{-w^2/2t} dt,$$

we have

$$I_\nu(z) = \int_0^\infty e^{-\nu^2 t/2} \bar{\theta}(r,t) dt,$$

$$\bar{\theta}(r,t) = \frac{1}{2\pi\sqrt{-1}} \int_C e^{r\cosh(w)} \frac{w}{\sqrt{2\pi t^3}} e^{-w^2/2t} dw.$$

We write down the complex integral as the sum of the integrals along the half lines and the line segment to obtain

$$\sqrt{2\pi t^3}\bar{\theta}(r,t)$$
$$= \frac{-1}{2\pi\sqrt{-1}} \int_0^\infty e^{-r\cosh(\xi)} e^{-(\xi-\pi\sqrt{-1})^2/2t} (\xi - \pi\sqrt{-1}) d\xi$$
$$- \frac{1}{2\pi\sqrt{-1}} \int_{-\pi}^\pi e^{r\cos(\xi)} e^{\xi^2/2t} \xi d\xi$$
$$+ \frac{1}{2\pi\sqrt{-1}} \int_0^\infty e^{-r\cosh(\xi)} e^{-(\xi+\pi\sqrt{-1})^2/2t} (\xi + \pi\sqrt{-1}) d\xi.$$

The second term on the right hand side is zero because the integrand is an odd function. Hence, we obtain after a simple computation

$$\sqrt{2\pi t^3}\bar{\theta}(r,t)$$
$$= \frac{1}{\pi} e^{\pi^2/2t} \int_0^\infty e^{-r\cosh(\xi)} e^{-\xi^2/2t} \left(\pi \cos\left(\frac{\pi\xi}{t}\right) - \xi \sin\left(\frac{\pi\xi}{t}\right)\right) d\xi,$$
$$= \frac{rt}{\pi} e^{\pi^2/2t} \int_0^\infty e^{-r\cosh(\xi)} \sinh(\xi) e^{-\xi^2/2t} \sin\left(\frac{\pi\xi}{t}\right) d\xi,$$

where we have used integration by parts formula for the second equality.

Now, by the injectivity of Laplace transform, it is easy to show $\bar{\theta}(r,t) = \theta(r,t)$. □